\documentclass[a4paper,11pt]{amsart}
\usepackage{graphicx}
\usepackage{mathptmx}
\usepackage{amsmath}
\usepackage{amssymb}
\usepackage{enumitem}
\usepackage{xcolor}

\newmuskip\pFqmuskip

\newcommand*\pFq[6][8]{%
  \begingroup 
  \pFqmuskip=#1mu\relax
  \mathcode`=\string"8000
  \begingroup\lccode`\~=`\,
  \lowercase{\endgroup\let~}\pFqcomma
  F^{#2}_{#3}{\left(\genfrac..{0pt}{}{#4}{#5}\bigg|#6\right)}%
  \endgroup
}
\newcommand{\pFqcomma}{\mskip\pFqmuskip}

\newtheorem{theorem}{Theorem}

\begin{document}

\title[Representations of degenerate Hermite polynomials]{Representations of degenerate Hermite polynomials }

\author{Taekyun  Kim}
\address{Department of Mathematics, Kwangwoon University, Seoul 139-701, Republic of Korea}
\email{tkkim@kw.ac.kr}

\author{Dae San Kim}
\address{Department of Mathematics, Sogang University, Seoul 121-742, Republic of Korea}
\email{dskim@sogang.ac.kr}

\author{Lee-Chae Jang}
\address{Graduate School of Eduction, Konkuk University, Seoul, 143-701, Republic of Korea}
\email{lcjang@konkuk.ac.kr}

\author{Hyunseok Lee}
\address{Department of Mathematics, Kwangwoon University, Seoul 139-701, Republic of Korea}
\email{luciasconstant@kw.ac.kr}

\author{Hanyoung Kim}
\address{Department of Mathematics, Kwangwoon University, Seoul 139-701, Republic of Korea}
\email{gksaud213@kw.ac.kr}


\maketitle

\begin{abstract}
We introduce degenerate Hermite polynomials as a degenerate version of the ordinary Hermite polynomials. Then, among other things, by using the formula about representing one $\lambda$-Sheffer polynomial in terms of other $\lambda$-Sheffer polynomials we represent the degenerate Hermite polynomials in terms of the higher-order degenerate Bernoulli, Euler, and Frobenius-Euler polynomials and vice versa.
\end{abstract}

\section{Introduction and preliminaries}
Carlitz initiated the study of degenerate Bernoulli and Euler polynomials and numbers and obtained some interesting arithmetical and combinatorial results in [3,4].
In recent years, some mathematicians began to investigate degenerate versions of quite a few special numbers and polynomials, which include the degenerate Bernoulli polynomials of the second kind, degenerate Stirling numbers of the first and second kinds, degenerate Bell polynomials, degenerate Frobenius-Euler polynomials, and so on (see [8,9,11-14] and the references therein).
It is remarkable that studying degenerate versions is not only limited to polynomials but also extended to transcendental functions. Indeed, the degenerate gamma functions were introduced in connection with degenerate Laplace transforms in [10]. \par
Gian-Carlo Rota began to construct a completely rigorous foundation for umbral calculus in the 1970s.
The Rota's theory is based on the linear functionals in \eqref{10} and \eqref{11} and differential operators in \eqref{14} and \eqref{15}. The Sheffer sequences, which are defined by \eqref{19}, occupy the central position in the theory and are characterized by the generating functions as in \eqref{20}, where the usual exponential function enters. The motivation for [8] started from the question that what if the usual exponential function in \eqref{20} was replaced by the degenerate exponential functions in \eqref{3-1}. It turns out that it corresponds to replacing the linear functional in \eqref{10} and \eqref{11} by the family of $\lambda$-linear functionals in \eqref{8} and \eqref{9} and the differential operators in \eqref{14} and \eqref{15} by the family of  $\lambda$-differential operators in \eqref{12} and \eqref{13}. Indeed, with these replacements we were led to define $\lambda$-Sheffer polynomials, which are charactered by the desired generating functions (see \eqref{18}). \par
The aim of this paper is to introduce the degenerate Hermite polynomials as a degenerate version of the ordinary Hermite polynomials and to study their properties. Specifically, by using the formula \eqref{23} about representing one $\lambda$-Sheffer polynomial by other $\lambda$-Sheffer polynomials we will represent the degenerate Hermite polynomial in terms of three other degenerate polynomials, namely the higher-order degenerate Bernoulli, Euler and Frobenius-Euler polynomials and vice versa.   \par
The outline of this paper is as follows. In Section 1, we will briefly go over very basics about umbral calculus including $\lambda$--linear functionals, $\lambda$--differential operators, $\lambda$--Sheffer sequences and the important formula in \eqref{23}. For further details on these, we let the reader refer to [8]. In addition, we  will recall the definitions for the Hermite polynomials, the degenerate Stirling numbers of the first and second kinds, higher-order degenerate Bernoulli polynomials, higher-order degenerate Euler polynomials and higher-order degenerate Frobenius--Euler polynomials. In Section 2, we will introduce the degenerate Hermite polynomials in a natural way and derive an explicit expression for them. Then we will represent the higher-order degenerate Bernoulli, Euler and Frobenius-Euler polynomials in terms of the degenerate Hermite polynomials and vice versa.

\vspace{0.1cm}

From now on, unless otherwise stated, $\lambda$ is arbitrary but is a fixed non-zero real number. It is well known that the Hermite polynomials, $H_{n}(x),\ (n\ge 0)$, are defined by the exponential generating function as
\begin{equation}
e^{2xt-t^{2}}=\sum_{n=0}^{\infty}H_{n}(x)\frac{t^{n}}{n!},\quad(\mathrm{see}\ [18]).\label{1}
\end{equation}
From \eqref{1}, we note that
\begin{equation}
H_{n}(x)=(-1)^{n}e^{x^{2}}\frac{d^{n}}{dx^{n}}e^{-x^{2}},\quad(n\ge 0).\label{2}	
\end{equation}
In particular,
\begin{displaymath}
	H_{n}(x)=\frac{n!}{2\pi i}\oint e^{-t^{2}+2xt}t^{-n-1}dt,\quad(n\ge 0),
\end{displaymath}
where the contour encloses the origin and is traversed in a  counterclockwise direction. \par
Carlitz introduced the degenerate Bernoulli polynomials of order $r (\in\mathbb{N})$ given by
\begin{equation}
\bigg(\frac{t}{e_{\lambda}(t)-1}\bigg)^{r}e_{\lambda}^{x}(t)=\sum_{n=0}^{\infty}\beta_{n,\lambda}^{(r)}(x)\frac{t^{n}}{n!},\quad(\mathrm{see}\ [3,4]).\label{3}
\end{equation}
Here $e_{\lambda}^{x}(t)$ are the degenerate exponential functions given by
\begin{equation}
e_{\lambda}^{x}(t)=(1+\lambda t)^{\frac{x}{\lambda}},\quad e_{\lambda}(t)=e_{\lambda}^{1}(t)=(1+\lambda t)^{\frac{1}{\lambda}},\quad(\mathrm{see}\ [10,14,15,16]). \label{3-1}
\end{equation}
For $x=0$, $\beta_{n,\lambda}^{(r)}=\beta_{n,\lambda}^{(r)}(0)$ are called the degenerate Bernoulli numbers of order $r$. \par
Note that $\displaystyle\lim_{\lambda\rightarrow 0}\beta_{n,\lambda}^{(r)}(x)=B_{n}^{(r)}(x)$, where $B_{n}^{(r)}(x)$ are the ordinary Bernoulli polynomials of order $r$, (see [1-20]). \par
Also, he considered the degenerate Euler polynomials of order $r$ given by
\begin{equation}
\bigg(\frac{2}{e_{\lambda}(t)+1}\bigg)^{r}e_{\lambda}^{x}(t)=\sum_{n=0}^{\infty}\mathcal{E}_{n,\lambda}^{(r)}(x)\frac{t^{n}}{n!},\quad(\mathrm{see}\ [3,4]),\label{4}
\end{equation}
where $\mathcal{E}_{n,\lambda}^{(r)}=\mathcal{E}_{n,\lambda}^{(r)}(0)$ are called the degenerate Euler numbers of order $r$. \par
The degenerate logarithmic functions are defined by Kim-Kim as
\begin{equation}
\log_{\lambda}(1+t)=\sum_{n=0}^{\infty}\lambda^{n-1}(1)_{n,1/\lambda}\frac{t^{n}}{n!}=\frac{1}{\lambda}\big((1+t)^{\lambda}-1\big),\quad(\mathrm{see}\ [9]),\label{5}	
\end{equation}
where $(x)_{0,\lambda}=1$, $(x)_{n,\lambda}=x(x-\lambda)\cdots\big(x-(n-1)\lambda\big),\,\,(n \ge 1)$. \par
Note that
\begin{displaymath}
\log_{\lambda}\big(e_{\lambda}(1+t)\big)=e_{\lambda}\big(\log_{\lambda}(1+t)\big)=1+t.
\end{displaymath}
The degenerate Stirling numbers of the first kind are defined by Kim-Kim as
\begin{equation}
(x)_{n}=\sum_{l=0}^{n}S_{1,\lambda}(n,l)(x)_{l,\lambda},\quad(n\ge 0),\quad(\mathrm{see}\ [9,11,12,13]),\label{6}
\end{equation}
where $(x)_{0}=1$, $(x)_{n}=x(x-1)\cdots(x-n+1)$, $(n\ge 1)$. \par
As an inversion formula of \eqref{6}, the degenerate Stirling numbers of the second kind are defined as
\begin{displaymath}
(x)_{n,\lambda}=\sum_{l=0}^{n}S_{2,\lambda}(n,l)(x)_{l},\quad(n\ge 0),\quad(\mathrm{see}\ [9,13,14,16]).
\end{displaymath}
In [14], the degenerate Frobenius-Euler polynomials of order $r$ are defined by
\begin{equation}
\bigg(\frac{1-u}{e_{\lambda}(t)-u}\bigg)^{r}e_{\lambda}^{x}(t)=\sum_{n=0}^{\infty}h_{n,\lambda}^{(r)}(x|u)\frac{t^{n}}{n!},\label{7}	
\end{equation}
where $u\in\mathbb{C}$ with $u\ne 1$, and $h_{n,\lambda}^{(r)}(u)=h_{n,\lambda}^{(r)}(0|u)$ are called the Frobenius-Euler numbers of order $r$. \par
Let $\mathbb{C}$ be the field of complex numbers. Let $\mathcal{F}$ be the algebra of formal power series in $t$ over $\mathbb{C}$ given by
\begin{displaymath}
\mathcal{F}=\bigg\{f(t)=\sum_{k=0}^{\infty}a_{k}\frac{t^{k}}{k!}\ \bigg|\ a_{k}\in\mathbb{C}\bigg\},
\end{displaymath}
and let $\mathbb{P}=\mathbb{C}[x]$ be the algebra of polynomials in $x$ over  $\mathbb{C}$. \par
For $f(t)\in\mathcal{F}$, with $\displaystyle f(t)=\sum_{k=0}^{\infty}a_{k}\frac{t^{k}}{k!}\displaystyle$, we define the $\lambda$-linear functional $\langle f(t)|\cdot \rangle_{\lambda}$ on $\mathbb{P}$ by
\begin{equation}
\langle f(t)|(x)_{k,\lambda}\rangle_{\lambda}=a_{k},\quad (k\ge 0),\quad (\mathrm{see}\ [8]).\label{8}
\end{equation}
In other words, the $\lambda$-linear functional $\langle f(t)|\cdot \rangle_{\lambda}$ is the unique linear functional on $\mathbb{P}$ mapping $(x)_{k,\lambda}$ onto $a_k$.
From \eqref{8}, we note that
\begin{equation}
\langle t^{k}|(x)_{n,\lambda}\rangle_{\lambda}=n!\delta_{n,k},\quad (n,k\ge 0),\quad (\mathrm{see}\ [8]),\label{9}
\end{equation}
where $\delta_{n,k}$ is the Kronceker's symbol. \par
If $\lambda=0$, then the symbol $\langle \cdot|\cdot \rangle_{0}$ is simply denoted by  $\langle \cdot|\cdot \rangle$. With this notation, \eqref{8} and \eqref{9} respectively become the linear functionals used in Rota's theory (see [18]):
\begin{equation}
\langle f(t)|x^k \rangle=a_{k},\quad (k\ge 0), \label{10}
\end{equation}
and
\begin{equation}
\langle t^{k}|x^n\rangle=n!\delta_{n,k},\quad (n,k\ge 0). \label{11}
\end{equation}
The order $o(f(t))$ of the formal power series $f(t)(\ne 0)$ is the smallest integer $k$ for which $a_{k}$ does not vanish. If $o(f(t))=1$, then $f(t)$ is called a delta series; if $o(f(t))=0$, then $f(t)$ is called an invertible series, (see [8,18]). \par
For each non-negative integer $k$, we define the $\lambda$-differential operator on $\mathbb{P}$ by
\begin{equation}
(t^{k})_{\lambda}(x)_{n,\lambda}=\left\{\begin{array}
	{cc}
	(n)_{k}(x)_{n-k,\lambda}, & \textrm{if $k\le n$,}\\
	0, &\textrm{if $k>n$,}
\end{array}\right.\quad(\mathrm{see}\ [8]),\label{12}
\end{equation}
and, extending this linearly, any power series
\begin{displaymath}
f(t)=\sum_{k=0}^{\infty}\frac{a_{k}}{k!}t^{k}\in\mathcal{F}
\end{displaymath}
yields the $\lambda$-differential operator on $\mathbb{P}$ given by
\begin{equation}
\big(f(t)\big)_{\lambda}(x)_{n,\lambda}=\sum_{k=0}^{n}\binom{n}{k}a_{k}(x)_{n-k,\lambda},\quad(n\ge 0).\label{13}
\end{equation}
In other words, we have
\begin{displaymath}
\big(f(t)\big)_{\lambda}=\sum_{k=0}^{\infty}\frac{a_{k}}{k!}(t^{k})_{\lambda}.
\end{displaymath} \par
If $\lambda=0$, then the differential operator $(t^k)_{0}$ is simply denoted by $t^k$. Now, \eqref{12} and \eqref{13} respectively become the differential operators used in Rota's theory:
\begin{equation}
t^{k}x^n=\left\{\begin{array}
{cc}
(n)_{k}x^{n-k}, & \textrm{if $k\le n$,}\\
0, &\textrm{if $k>n$,}
\end{array}\right. \quad(\mathrm{see}\ [18]),\label{14}
\end{equation}
and, extending this linearly, any formal power series
\begin{equation*}
f(t)=\sum_{k=0}^{\infty}\frac{a_{k}}{k!}t^{k}\in\mathcal{F}
\end{equation*}
yields the differential operator on $\mathbb{P}$ given by
\begin{equation}
f(t)x^n=\sum_{k=0}^{n}\binom{n}{k}a_{k}x^{n-k},\quad(n\ge 0).\label{15}
\end{equation}
From \eqref{13}, we note that
\begin{equation}
\big(e_{\lambda}^{y}(t)\big)_{\lambda}(x)_{n,\lambda}=(x+y)_{n,\lambda},\quad\big(e_{\lambda}^{y}(t)\big)_{\lambda}P(x)=P(x+y),\quad P(x)\in\mathbb{P},\label{16}
\end{equation}
\begin{displaymath}
\big(e_{\lambda}^{y}(t)-1\big)_{\lambda}P(x)=P(x+y)-P(x),\quad\big\langle e_{\lambda}^{y}(t)|P(x)\rangle_{\lambda}=P(y),\quad (\mathrm{see}\ [8]).
\end{displaymath}
For $f(t),g(t)\in\mathcal{F}$ with $o(f(t))=1$ and $o(g(t))=0$, there exists a unique sequences $s_{n,\lambda}(x)$ $(\deg s_{n,\lambda}(x)=n)$ of polynomials such that (see [8])
\begin{equation}
\big\langle g(t)\big(f(t)\big)^{k}|s_{n,\lambda}(x)\big\rangle_{\lambda}=n!\delta_{n,k},\quad(n,k\ge 0). \label{17}
\end{equation} \par
The sequence $s_{n,\lambda}(x)$ is called the $\lambda$-Sheffer sequence for $(g(t),f(t))$, which is denoted by $s_{n,\lambda}(x)\sim (g(t),f(t))_{\lambda}$, (see [8]). \par
For $s_{n,\lambda}(x)\sim (g(t),f(t))_{\lambda}$, we have
\begin{equation}
\frac{1}{g(\overline{f}(t))}e_{\lambda}^{y}\big(\overline{f}(t)\big)=\sum_{k=0}^{\infty}s_{k,\lambda}(y)\frac{t^{k}}{k!},\quad(\mathrm{see}\ [8]),\label{18}
\end{equation}
for all $y\in\mathbb{C}$, where $\overline{f}(t)$ is the compositional inverse of $f(t)$ with $\overline{f}(f(t))=f(\overline{f}(t))=t$. \par
For $\lambda=0$, and with $f(t),g(t)\in\mathcal{F}$ as before, there exists a unique sequence $s_n(x)$  $(\deg s_{n}(x)=n)$ of polynomials such that  (see [18])
\begin{equation}
\big\langle g(t)\big(f(t)\big)^{k}|s_{n}(x)\big\rangle=n!\delta_{n,k},\quad(n,k\ge 0). \label{19}
\end{equation} \par
The sequence $s_{n}(x)$ is called the Sheffer sequence for $(g(t),f(t))$, which is denoted by $s_{n}(x)\sim (g(t),f(t))$, (see [18]). \par
For $s_{n}(x)\sim (g(t),f(t))$, we have
\begin{equation}
\frac{1}{g(\overline{f}(t))}e^{y \overline{f}(t)}=\sum_{k=0}^{\infty}s_{k}(y)\frac{t^{k}}{k!},\quad(\mathrm{see}\ [18]), \label{20}
\end{equation}
for all $y\in\mathbb{C}$, where $\overline{f}(t)$ is the compositional inverse of $f(t)$ with $\overline{f}(f(t))=f(\overline{f}(t))=t$. \par
Assume that, for each $ \lambda \in \mathbb{R}^{*}$ of the set of nonzero real numbers, $s_{n,\lambda}(x)$ is $\lambda$-Sheffer for $(g_{\lambda}(t), f_{\lambda}(t))$. Assume also that $\lim_{\lambda \rightarrow 0} f_{\lambda}(t)=f(t),\,\, \lim_{\lambda \rightarrow 0} g_{\lambda}(t)=g(t)$, for some delta series $f(t)$ and an invertible series $g(t)$. Then we see that $\lim_{\lambda \rightarrow 0} \overline{f}_{\lambda}(t)=\overline{f}(t)$. Moreover, by \eqref{18}, for each  $ \lambda \in \mathbb{R}^{*}$ we have:
\begin{equation}
\frac{1}{g_{\lambda}\big(\overline{f}_{\lambda}(t)\big)}e_{\lambda}^{x}\big(\overline{f}_{\lambda}(t)\big)\ =\ \sum_{k=0}^{\infty}s_{k,\lambda}(x)\frac{t^{k}}{k!}.\label{21}
\end{equation}
If $ \lim_{\lambda \rightarrow 0}s_{k,\lambda}(x)=s_k(x)$, then, by \eqref{21}, we have
\begin{equation*}
\frac{1}{g\big(\overline{f}(t)\big)}e^{x\overline{f}(t)}\ =\ \sum_{k=0}^{\infty}s_{k}(x)\frac{t^{k}}{k!}.
\end{equation*}
Hence $s_n(x)$ is Sheffer for $(g(t), f(t))$, and
\begin{equation*}
\big\langle g(t)f(t)^{k}|s_{n}(x)\big\rangle=n!\delta_{n,k},\quad(n,k\ge 0).
\end{equation*}
In this case, we may say that the family $\left\{s_{n,\lambda}(x) \right\}_{\lambda \in \mathbb{R}^{*}}$ of $\lambda$-Sheffer sequences $s_{n,\lambda}(x)$ are the degenerate sequences for the Sheffer polynomial $s_n(x)$. For example, $\left\{\beta_{n,\lambda}(x) \right\}_{\lambda \in \mathbb{R}^{*}}$ (see \eqref{3}, with $r=1$) are the degenerate sequences for the Bernoulli polynomials $B_n(x)$, where $B_n(x)=\lim_{\lambda \rightarrow 0}\beta_{n,\lambda}(x)$ are given by
\begin{equation*}
\frac{t}{e^t -1} e^{xt}=\sum_{k=0}^{\infty}B_k(x) \frac{t^k}{k!}.
\end{equation*}
However, $s_{n,\lambda}(x)$ itself is oftentimes called the degenerate Sheffer polynomial for the Sheffer polynomial $s_n(x)$. \par
For $f(t),g(t)\in\mathcal{F}$ and $P(x)\in\mathbb{P}$, it is easy to show that
\begin{displaymath}
\langle f(t)g(t)|P(x)\rangle_{\lambda}=\langle g(t)|\big(f(t)\big)_{\lambda}P(x)\rangle_{\lambda}=\langle f(t)|\big(g(t)\big)_{\lambda}P(x)\rangle_{\lambda},\quad(\mathrm{see}\ [8]).
\end{displaymath}
For $s_{n,\lambda}(x)\sim (g(t),f(t))_{\lambda}$, $r_{n,\lambda}(x)\sim (h(t),l(t))_{\lambda}$, we have
\begin{equation}
s_{n,\lambda}(x)=\sum_{k=0}^{n}c_{n,k}r_{k,\lambda}(x),\quad(\mathrm{see}\ [8]), \label{22}
\end{equation}
where
\begin{equation}
c_{n,k}=\frac{1}{k!}\bigg\langle\frac{h(\overline{f}(t))}{g(\overline{f}(t))}\big(l(\overline{f}(t))\big)^{k}\bigg|(x)_{n,\lambda}\bigg\rangle_{\lambda}. \label{23}
\end{equation}

\section{Representations of degenerate Hermite polynomials in terms of other degenerate Sheffer polynomials and vice versa}
In light of \eqref{1}, we may consider the degenerate Hermite polynomials which are given by
\begin{equation}
e_{\lambda}^{-1}(t^{2})\cdot e_{\lambda}^{x}(2t)=\sum_{n=0}^{\infty}H_{n,\lambda}(x)\frac{t^{n}}{n!}. \label{24}
\end{equation}
Note that
\begin{displaymath}
\sum_{n=0}^{\infty}\lim_{\lambda\rightarrow 0}H_{n,\lambda}(x)\frac{t^{n}}{n!}=e^{-t^{2}+2xt}=\sum_{n=0}^{\infty}H_{n}(x)\frac{t^{n}}{n!}.
\end{displaymath}
By \eqref{24}, we get
\begin{align}
\sum_{n=0}^{\infty}H_{n,\lambda}(x)\frac{t^{n}}{n!}\ &=\ e_{\lambda}^{-1}(t^{2})\cdot e_{\lambda}^{x}(2t)\ =\ \sum_{l=0}^{\infty}(-1)_{l,\lambda}\frac{t^{2l}}{l!}\sum_{m=0}^{\infty}(x)_{m,\lambda}2^m\frac{t^{m}}{m!} \label{25} \\
&=\ \sum_{n=0}^{\infty}\bigg(n!\sum_{l=0}^{\big[\frac{n}{2}\big]}\frac{(-1)_{l,\lambda}2^{n-2l}(x)_{n-2l,\lambda}}{l!(n-2l)!}\bigg)\frac{t^{n}}{n!}.\nonumber
\end{align}
By comparing the coefficients on both sides of \eqref{25}, we obtain the following theorem.
\begin{theorem}
For $n\ge 0$, we have the expression given by
\begin{displaymath}
H_{n,\lambda}(x)=n!\sum_{l=0}^{\big[\frac{n}{2}\big]}\frac{2^{n-2l} (-1)_{l,\lambda}}{l!(n-2l)!} (x)_{n-2l,\lambda},
\end{displaymath}
where $[x]$ denotes the greatest integer not exceeding $x$.
\end{theorem}
Note that $H_{n,\lambda}(x)$ is the $\lambda$-Sheffer sequence for $\big(e_{\lambda}\big(\frac{1}{4}t^{2}\big),\frac{t}{2}\big)$. Now, we consider the following two $\lambda$-Sheffer sequences:
\begin{equation}
\mathcal{E}_{n,\lambda}^{(r)}(x)\sim\bigg(\bigg(\frac{e_{\lambda}(t)+1}{2}\bigg)^{r}, t\bigg)_{\lambda},\quad H_{n,\lambda}(x)\sim\bigg(e_{\lambda}\bigg(\frac{1}{4}t^{2}\bigg),\frac{1}{2}t\bigg)_{\lambda}.\label{26}
\end{equation}
By \eqref{22}, \eqref{23} and \eqref{26}, we get
\begin{equation}
\mathcal{E}_{n,\lambda}^{(r)}(x)=\sum_{k=0}^{n}C_{n,k}H_{n,\lambda}(x),\label{27}
\end{equation}
where
\begin{align}
C_{n,k}\ &=\ \frac{1}{k!}\bigg\langle\frac{e_{\lambda}\big(\frac{1}{4}t^{2}\big)}{\big(\frac{e_{\lambda}(t)+1}{2}\big)^{r}}\bigg(\frac{t}{2}\bigg)^{k}\bigg|(x)_{n,\lambda}\bigg\rangle_{\lambda}\label{28} \\
&=\ \frac{1}{2^{k}k!}\bigg\langle\bigg(\frac{2}{e_{\lambda}(t)+1}\bigg)^r e_{\lambda}\bigg(\frac{1}{4}t^{2}\bigg)\bigg|(t^{k})_{\lambda}(x)_{n,\lambda}\bigg\rangle_{\lambda}\nonumber 	\\
&=\ \frac{\binom{n}{k}}{2^{k}}\bigg\langle\bigg(\frac{2}{e_{\lambda}(t)+1}\bigg)^r \bigg|\bigg(e_{\lambda}\bigg(\frac{1}{4}t^{2}\bigg)\bigg)_{\lambda}(x)_{n-k,\lambda}\bigg\rangle_{\lambda}\nonumber\\
&=\ 2^{-k}\binom{n}{k}\sum_{l=0}^{\big[\frac{n-k}{2}\big]}\frac{(1)_{l,\lambda}}{4^{l}l!}(n-k)_{2l} \bigg\langle\bigg(\frac{2}{e_{\lambda}(t)+1}\bigg)^{r}\bigg|(x)_{n-k-2l,\lambda}\bigg\rangle_{\lambda}\nonumber\\
&=\ 2^{-k}\binom{n}{k}\sum_{l=0}^{\big[\frac{n-k}{2}\big]}\frac{(1)_{l,\lambda}}{2^{2l}l!}(n-k)_{2l}\mathcal{E}_{n-k-2l,\lambda}^{(r)} \nonumber\\
&=\ n!\sum_{\substack{0\le l\le n-k\\ l:\ even}}\frac{(1)_{\frac{l}{2},\lambda}}{2^{l+k}\big(\frac{l}{2}\big)!(n-k-l)!k!}\mathcal{E}_{n-k-l,\lambda}^{(r)}.\nonumber
\end{align}

Therefore, by \eqref{27} and \eqref{28}, we obtain the following theorem.
\begin{theorem}
For $n\ge 0$, we have the representation given by
\begin{displaymath}
\mathcal{E}_{n,\lambda}^{(r)}(x)=n!\sum_{k=0}^{n}\bigg\{\sum_{\substack{0\le l\le n-k\\ l:\ even}}\frac{(1)_{\frac{l}{2},\lambda}}{k!(n-k-l)!2^{l+k}\big(\frac{l}{2}\big)!}\mathcal{E}_{n-k-l,\lambda}^{(r)}\bigg\}H_{n,\lambda}(x).
\end{displaymath}
\end{theorem}
\noindent For
\begin{equation}
\beta_{n,\lambda}^{(r)}(x)\sim\bigg(\bigg(\frac{e_{\lambda}(t)-1}{t}\bigg)^{r},t\bigg)_{\lambda},\quad H_{n,\lambda}(x)\sim\bigg(e_{\lambda}\bigg(\frac{1}{4}t^{2}\bigg),\frac{t}{2}\bigg)_{\lambda},\label{29}
\end{equation}
we have
\begin{equation}
\beta_{n,\lambda}^{(r)}(x)=\sum_{k=0}^{n}C_{n,k}H_{k,\lambda}(x),\label{30}
\end{equation}
where
\begin{align}
C_{n,k}\ &=\ \frac{1}{k!}\bigg\langle\frac{e_{\lambda}\big(\frac{1}{4}t^{2}\big)}{\big(\frac{e_{\lambda}(t)-1}{t}\big)^{r}}\bigg(\frac{t}{2}\bigg)^{k}\bigg|(x)_{n,\lambda}\bigg\rangle_{\lambda}\label{31} \\
&=\ \frac{1}{k!2^{k}}\bigg\langle \bigg(\frac{t}{e_{\lambda}(t)-1}\bigg)^{r}e_{\lambda}\bigg(\frac{1}{4}t^{2}\bigg)\bigg|(t^{k})_{\lambda}(x)_{n,\lambda}\bigg\rangle_{\lambda}\nonumber \\
&=\ 2^{-k}\binom{n}{k}\sum_{l=0}^{\big[\frac{n-k}{2}\big]}\frac{(1)_{l,\lambda}}{l!4^{l}}(n-k)_{2l}\bigg\langle\bigg(\frac{t}{e_{\lambda}(t)-1}\bigg)^{r}\bigg|(x)_{n-k-2l,\lambda}\bigg\rangle_{\lambda}\nonumber
\end{align}
\begin{align*}
&=\ 2^{-k}\binom{n}{k}\sum_{l=0}^{\big[\frac{n-k}{2}\big]}\frac{(1)_{l,\lambda}}{l!4^{l}}(n-k)_{2l}\beta_{n-k-2l,\lambda}^{(r)}\nonumber\\
&=\ n!\sum_{\substack{0\le l\le n-k\\ l:\ even}}\frac{(1)_{\frac{l}{2},\lambda}}{k!(n-k-l)!2^{k+l}\big(\frac{l}{2}\big)!}\beta_{n-k-l,\lambda}^{(r)}.\nonumber
\end{align*}
Therefore, by \eqref{30} and \eqref{31}, we obtain the following theorem.
\begin{theorem}
For $n\ge 0$, we have the representation given by
\begin{displaymath}
\beta_{n,\lambda}^{(r)}(x)=n!\sum_{k=0}^{n}\bigg\{\sum_{\substack{0\le l\le n-k\\ l:\ even}}\frac{(1)_{\frac{l}{2},\lambda} \beta_{n-k-l,\lambda}^{(r)}}{k!(n-k-l)!\big(\frac{l}{2}\big)! 2^{k+l}} \bigg\}H_{k,\lambda}(x).
\end{displaymath}
\end{theorem}
From \eqref{7}, \eqref{18} and \eqref{24}, we note that
\begin{equation}
h_{n,\lambda}^{(r)}(x|u)\sim\bigg(\bigg(\frac{e_{\lambda}(t)-u}{1-u}\bigg)^{r},t\bigg)_{\lambda},\quad H_{n,\lambda}(x)\sim \bigg(e_{\lambda}\bigg(\frac{1}{4}t^{2}\bigg),\frac{t}{2}\bigg)_{\lambda}.\label{32}	
\end{equation}
By \eqref{22}, \eqref{23} and \eqref{32}, we get
\begin{equation}
h_{n,\lambda}^{(r)}(x|u)=\sum_{k=0}^{n}C_{n,k}H_{k,\lambda}(x),\label{33}	
\end{equation}
where
\begin{align}
C_{n,k}\ &=\ \frac{1}{k!}\bigg\langle\frac{e_{\lambda}\big(\frac{1}{4}t^{2}\big)}{\big(\frac{e_{\lambda}(t)-u}{1-u}\big)^{r}}\bigg(\frac{t}{2}\bigg)^{k}\bigg|(x)_{n,\lambda}\bigg\rangle_{\lambda}\label{34} \\
&=\ \frac{1}{2^{k}k!}\bigg\langle \bigg(\frac{1-u}{e_{\lambda}(t)-u}\bigg)^{r}e_{\lambda}\bigg(\frac{1}{4}t^{2}\bigg)\bigg|(t^{k})_{\lambda}(x)_{n,\lambda}\bigg\rangle_{\lambda}\nonumber\\
&=\ \frac{\binom{n}{k}}{2^{k}} \bigg\langle \bigg(\frac{1-u}{e_{\lambda}(t)-u}\bigg)^{r}\bigg|\bigg(e_{\lambda}\bigg(\frac{1}{4}t^{2}\bigg)\bigg)_{\lambda}(x)_{n-k,\lambda}\bigg\rangle_{\lambda}\nonumber \\
&=\ \frac{\binom{n}{k}}{2^{k}}\sum_{l=0}^{\big[\frac{n-k}{2}\big]}\frac{(n-k)_{2l}(1)_{l,\lambda}}{l!2^{2l}}\bigg\langle\bigg(\frac{1-u}{e_{\lambda}(t)-u}\bigg)^{r}\bigg|(x)_{n-k-2l,\lambda}\bigg\rangle_{\lambda}\nonumber \\
&=\ \frac{\binom{n}{k}}{2^{k}}\sum_{l=0}^{\big[\frac{n-k}{2}\big]} \frac{(n-k)_{2l}(1)_{l,\lambda}}{l!2^{2l}}h_{n-k-2l,\lambda}^{(r)}(u)\nonumber \\
&=\ n!\sum_{\substack{0\le l\le n-k\\ l;\ even}}\frac{(1)_{\frac{l}{2},\lambda}h_{n-k-l,\lambda}^{(r)}(u)}{\big(\frac{l}{2}\big)!2^{l+k}(n-k-l)!k!}.\nonumber
\end{align}
Therefore, by \eqref{33} and \eqref{34}, we obtain the following theorem.
\begin{theorem}
For $n\ge 0$, we have the representation given by
\begin{displaymath}
h_{n,\lambda}^{(r)}(x|u)=n!\sum_{k=0}^{n}\bigg\{\ \sum_{\substack{0\le l\le n-k\\ l;\ even}}\frac{(1)_{\frac{l}{2},\lambda}h_{n-k-l,\lambda}^{(r)}(u)}{\big(\frac{l}{2}\big)!2^{l+k}(n-k-l)!k!}\bigg\}H_{k,\lambda}(x).
\end{displaymath}
\end{theorem}
For $s_{n,\lambda}(x)\sim (g(t),f(t))_{\lambda}$, $p_{n,\lambda}(x)\sim (1,f(t))_{\lambda}$, we have
\begin{align}
n!\delta_{n,k}\ &=\ \big\langle g(t)\big(f(t)\big)^{k}\big|s_{n,\lambda}(x)\big\rangle_{\lambda}\ =\ \big\langle \big(f(t)\big)^{k}\big|\big(g(t)\big)_{\lambda}s_{n,\lambda}(x)\big\rangle_{\lambda} \nonumber \\
&=\ \big\langle \big(f(t)\big)^{k}\big|p_{n,\lambda}(x)\big\rangle_{\lambda}\label{35}	
\end{align}
By \eqref{35}, we get
\begin{equation}
\big(g(t)\big)_{\lambda}s_{n,\lambda}(x)\sim (1,f(t))_{\lambda},\quad (n\ge 0). \label{36}
\end{equation}
From \eqref{24}, we have
\begin{equation}
H_{n,\lambda}(x)\sim\bigg(e_{\lambda}\bigg(\frac{1}{4}t^{2}\bigg),\frac{1}{2}t\bigg)_{\lambda}\Longleftrightarrow\bigg(e_{\lambda}\bigg(\frac{1}{4}t^{2}\bigg)\bigg)_{\lambda}H_{n,\lambda}(x) \sim\bigg(1,\frac{1}{2}t\bigg)_{\lambda}  \label{37}.	
\end{equation}
Since $2^{n}(x)_{n,\lambda}\sim\big(1,\frac{1}{2}t)_{\lambda}$, $\big(e_{\lambda}\big(\frac{1}{4}t^{2}\big)\big)_{\lambda}H_{n,\lambda}(x)=2^{n}(x)_{n,\lambda}$. \par
Hence,
\begin{align*}
	H_{n,\lambda}(x)\ &=\ 2^{n}\bigg(e_{\lambda}^{-1}\bigg(\frac{1}{4}t^{2}\bigg)\bigg)_{\lambda}(x)_{n,\lambda}\\
	&=\ 2^{n}\sum_{l=0}^{\infty}\frac{(-1)_{l,\lambda}}{l!4^{l}}\big(t^{2l}\big)_{\lambda}(x)_{n,\lambda}\\
	&=\ 2^{n}\sum_{l=0}^{[\frac{n}{2}]}\frac{(-1)_{l,\lambda}(n)_{2l}}{l!4^{l}}(x)_{n-2l,\lambda}\\
	&=\ n!\sum_{l=0}^{[\frac{n}{2}]}\frac{(-1)_{l,\lambda}2^{n-2l}}{l!(n-2l)!}(x)_{n-2l,\lambda},
\end{align*}
which gives another proof for Theorem 1. Assume that
\begin{equation}
H_{n,\lambda}(x)=\sum_{k=0}^{n}C_{n,k}\mathcal{E}_{k,\lambda}^{(r)}(x).\label{38}	
\end{equation}
Then, by \eqref{22}, \eqref{23} and \eqref{26}, we get
\begin{align}
C_{n,k}\ &=\ \frac{1}{k!}\bigg\langle\frac{\big(\frac{e_{\lambda}(2t)+1}{2}\big)^{r}}{e_{\lambda}\big(\frac{1}{4}(2t)^{2}\big)}(2t)^{k}\bigg|(x)_{n,\lambda}\bigg\rangle_{\lambda}\label{39} \\
&=\frac{2^{k}}{k!}\bigg\langle\bigg(\frac{e_{\lambda}(2t)+1}{2}\bigg)^{r}e_{\lambda}^{-1}(t^{2})\bigg|(t^{k})_{\lambda}(x)_{n,\lambda}\bigg\rangle_{\lambda}\nonumber \\
&=\ 2^{k}\binom{n}{k}\bigg\langle\bigg(\frac{e_{\lambda}(2t)+1}{2}\bigg)^{r}\bigg|\big(e_{\lambda}^{-1}(t^{2})\big)_{\lambda}(x)_{n-k,\lambda}\bigg\rangle_{\lambda}\nonumber \\
&=\ 2^{k}\binom{n}{k}\sum_{l=0}^{[\frac{n-k}{2}]}\frac{(-1)_{l,\lambda}(n-k)_{2l}}{l!}\bigg\langle \bigg(\frac{e_{\lambda}(2t)+1}{2}\bigg)^{r}\bigg|(x)_{n-k-2l,\lambda}\bigg\rangle_{\lambda}\nonumber
\end{align}
\begin{align*}
&=\ \frac{2^{k}}{2^{r}}\binom{n}{k}\sum_{j=0}^{r}\sum_{l=0}^{[\frac{n-k}{2}]}\frac{(-1)_{l,\lambda}(n-k)_{2l}}{l!}\binom{r}{j}\big\langle e^{j}_{\lambda}(2t)\big|(x)_{n-k-2l,\lambda}\big\rangle_{\lambda}\nonumber\\
&=\ \frac{1}{2^{r}}\sum_{j=0}^{r}\sum_{l=0}^{[\frac{n-k}{2}]}\frac{(-1)_{l,\lambda}(n-k)_{2l}\binom{r}{j}\binom{n}{k}2^{k}2^{n-k-2l}(j)_{n-k-2l,\lambda}}{l!} \nonumber\\
&=\ \frac{1}{2^{r}}\sum_{j=0}^{r}\sum_{l=0}^{[\frac{n-k}{2}]}\frac{(-1)_{l,\lambda}(n-k)!\binom{r}{j}\binom{n}{k}2^{n-2l}}{l!(n-k-2l)!} (j)_{n-k-2l,\lambda}.\nonumber
\end{align*}
Therefore, by \eqref{38} and \eqref{39}, we obtain the following theorem.
\begin{theorem}
For $n\ge 0$, we have the representation given by
\begin{displaymath}
H_{n,\lambda}(x)=\frac{1}{2^{r}}\sum_{k=0}^{n}\bigg\{\sum_{j=0}^{r}\sum_{l=0}^{[\frac{n-k}{2}]}\frac{(-1)_{l,\lambda}(n-k)!\binom{r}{j}\binom{n}{k}2^{n-2l}}{l!(n-k-2l)!} (j)_{n-k-2l,\lambda}\bigg\}\mathcal{E}_{k,\lambda}^{(r)}(x).
\end{displaymath}
\end{theorem}
Next, we would like to find another representation of $H_{n,\lambda}(x)$ in terms of $\mathcal{E}_{k,\lambda}^{(r)}(x)$.
From \eqref{39}, we observe that
\begin{align}
C_{n,k}\ &=\ \frac{1}{k!}\bigg\langle\bigg(\frac{e_{\lambda}(2t)+1}{2}\bigg)^{r} e_{\lambda}^{-1}\bigg(\frac{1}{4}(2t)^{2}\bigg)(2t)^{k}\bigg|(x)_{n,\lambda}\bigg\rangle_{\lambda}\label{40}\\
&=\ \frac{1}{k!}\bigg\langle \bigg(\frac{e_{\lambda}(t)+1}{2}\bigg)^{r}e_{\lambda}^{-1}\bigg(\frac{1}{4}t^{2}\bigg)t^{k}\bigg|2^{n}(x)_{n,\lambda}\bigg\rangle_{\lambda}\nonumber\\
&=\ \frac{1}{k!}\bigg\langle \bigg(\frac{e_{\lambda}(t)+1}{2}\bigg)^{r}t^{k}\bigg|2^{n} \bigg(e_{\lambda}^{-1}\bigg(\frac{1}{4}t^{2}\bigg)\bigg)_{\lambda}(x)_{n,\lambda}\bigg\rangle_{\lambda}\nonumber \\
&=\ \frac{1}{k!}\bigg\langle \bigg(\frac{e_{\lambda}(t)+1}{2}\bigg)^{r}t^{k}\bigg|H_{n,\lambda}(x)\bigg\rangle_{\lambda}\nonumber \\
&=\ \frac{1}{k!}\bigg\langle\bigg(\frac{e_{\lambda}(t)+1}{2}\bigg)^{r}\bigg|(t^{k})_{\lambda}H_{n,\lambda}(x)\bigg\rangle_{\lambda}.\nonumber
\end{align}
To proceed further, we note that
\begin{align}
(t^{k})_{\lambda}H_{n,\lambda}(x)\ &=\ n!\sum_{l=0}^{[\frac{n-k}{2}]}\frac{2^{n-2l}(-1)_{l,\lambda}}{l!(n-2l)!}(t^{k})_{\lambda}(x)_{n-2l,\lambda}\nonumber \\
&=\ n!\sum_{l=0}^{[\frac{n-k}{2}]}\frac{2^{n-2l}(-1)_{l,\lambda}}{l!(n-2l-k)!}(x)_{n-2l-k,\lambda} \label{41} \\
&=\ 2^{k}(n)_{k}(n-k)! \sum_{l=0}^{[\frac{n-k}{2}]}\frac{2^{n-k-2l}(-1)_{l,\lambda}}{l!(n-k-2l)!}(x)_{n-2l-k,\lambda}\nonumber \\
&=\ 2^{k}(n)_{k}H_{n-k,\lambda}(x). \nonumber
\end{align}
By \eqref{40} and \eqref{41}, we get
\begin{align}
C_{n,k}\ &=\ \frac{2^{k}(n)_{k}}{k!}\bigg\langle\bigg(\frac{e_{\lambda}(t)+1}{2}\bigg)^{r}\bigg|H_{n-k,\lambda}(x)\bigg\rangle_{\lambda}\label{42} \\
&=\frac{2^{k}\binom{n}{k}}{2^{r}}\sum_{j=0}^{r}\binom{r}{j}\big\langle e_{\lambda}^{j}(t)\big|H_{n-k,\lambda}(x)\big\rangle_{\lambda} \nonumber \\
&=\ \frac{2^{k}}{2^{r}}\binom{n}{k}\sum_{j=0}^{r}\binom{r}{j}H_{n-k,\lambda}(j). \nonumber
\end{align}
Therefore, by \eqref{38} and \eqref{42}, we obtain the following theorem.
\begin{theorem}
For $n\ge 0$, we have the representation given by
\begin{displaymath}
H_{n,\lambda}(x)=\frac{1}{2^{r}}\sum_{k=0}^{n}\binom{n}{k}2^{k}\bigg[\sum_{j=0}^{r}\binom{r}{j}H_{n-k,\lambda}(j)\bigg]\mathcal{E}_{k,\lambda}^{(r)}(x).
\end{displaymath}
\end{theorem}
Let us assume that
\begin{equation}
H_{n,\lambda}(x)=\sum_{k=0}^{n}C_{n,k}\beta_{n,\lambda}^{(r)}(x). \label{43}
\end{equation}
Then, by \eqref{22}, \eqref{23} and \eqref{29}, we have
\begin{align}
C_{n,k}\ &=\ \frac{1}{k!}\bigg\langle\frac{\big(\frac{e_{\lambda}(2t)-1}{2t}\big)^{r}}{e_{\lambda}\big(\frac{1}{4}(2t)^{2}\big)}(2t)^{k}\bigg|(x)_{n,\lambda}\bigg\rangle_{\lambda}\label{44} \\
&=\ \frac{1}{k!}\bigg\langle\frac{\big(\frac{e_{\lambda}(t)-1}{t}\big)^{r}}{e_{\lambda}\big(\frac{1}{4}t^{2}\big)}t^{k}\bigg|2^{n}(x)_{n,\lambda}\bigg\rangle_{\lambda}\nonumber \\
&=\ \frac{1}{k!}\bigg\langle\bigg(\frac{e_{\lambda}(t)-1}{t}\bigg)^{r}t^{k}\bigg|\bigg(e_{\lambda}^{-1}\bigg(\frac{1}{4}t^{2}\bigg)\bigg)_{\lambda}2^{n}(x)_{n,\lambda}\bigg\rangle_{\lambda}\nonumber \\
&=\frac{1}{k!} \bigg\langle\bigg(\frac{e_{\lambda}(t)-1}{t}\bigg)^{r}t^{k}\bigg|H_{n,\lambda}(x)\bigg\rangle_{\lambda}. \nonumber
\end{align}
For $r>k$, we have
\begin{align}
	C_{n,k}\ &=\ \frac{1}{k!}\bigg\langle\big(e_{\lambda}(t)-1\big)^{k}\bigg|\bigg(\bigg(\frac{e_{\lambda}(t)-1}{t}\bigg)^{r-k}\bigg)_{\lambda}H_{n,\lambda}(x)\bigg\rangle_{\lambda} \label{45}\\
	&=\ \frac{1}{k!}\sum_{l=0}^{n}\frac{S_{2,\lambda}(l+r-k,r-k)}{\binom{l+r-k}{l}l!}\big\langle\big(e_{\lambda}(t)-1\big)^{k}\big|(t^{l})_{\lambda}H_{n,\lambda}(x)\big\rangle_{\lambda}\nonumber \\
		&=\ \frac{1}{k!}\sum_{l=0}^{n}\frac{S_{2,\lambda}(l+r-k,r-k)}{\binom{l+r-k}{l}l!}\big\langle\big(e_{\lambda}(t)-1\big)^{k}\big|2^{l}(n)_{l}H_{n-l,\lambda}(x)\big\rangle_{\lambda}\nonumber \\
		&=\ \frac{1}{k!}\sum_{l=0}^{n}\sum_{j=0}^{k}\frac{2^{l}\binom{n}{l}\binom{k}{j}}{\binom{l+r-k}{l}}S_{2,\lambda}(l+r-k,r-k)(-1)^{k-j}\big\langle e_{\lambda}^{j}(t)\big|H_{n-l,\lambda}(x)\big\rangle_{\lambda}\nonumber \\
		&=\ \frac{1}{k!}\sum_{l=0}^{n}\sum_{j=0}^{k}\frac{2^{l}\binom{n}{l}\binom{k}{j}}{\binom{l+r-k}{l}}S_{2,\lambda}(l+r-k,r-k)(-1)^{k-j}H_{n-l,\lambda}(j).\nonumber
\end{align}
The next theorem now follows from \eqref{43} and \eqref{45}.
\begin{theorem}
	For $r>n\ge 0$, we have the representation given by
	\begin{displaymath}
		H_{n,\lambda}(x)=\sum_{k=0}^{n}\frac{1}{k!}\bigg\{\sum_{j=0}^{k}\sum_{l=0}^{n}\frac{2^{l}\binom{n}{l}\binom{k}{j}}{\binom{l+r-k}{l}}S_{2,\lambda}(l+r-k,r-k)(-1)^{k-j}H_{n-l,\lambda}(j)\bigg\}\beta_{k,\lambda}^{(r)}(x).
	\end{displaymath}
\end{theorem}
For $r\le k\le n$ in \eqref{44}, we have
\begin{align}
C_{n,k}\ &=\ \frac{1}{k!}\sum_{j=0}^{r}\binom{r}{j}(-1)^{r-j}\big\langle e_{\lambda}^{j}(t)t^{k-r}\big|H_{n,\lambda}(x)\rangle_{\lambda}\label{46} \\
&=\ \frac{1}{k!}\sum_{j=0}^{r}\binom{r}{j}(-1)^{r-j}\big\langle e_{\lambda}^{j}(t)\big|\big(t^{k-r}\big)_{\lambda}H_{n,\lambda}(x)\big\rangle_{\lambda}\nonumber \\
&=\ \frac{2^{k-r}(n)_{k-r}}{k!}\sum_{j=0}^{r}\binom{r}{j}(-1)^{r-j}\big\langle e_{\lambda}^{j}(t)\big|H_{n-k+r,\lambda}(x)\big\rangle_{\lambda}\nonumber\\
&=\ \frac{2^{k-r}(n)_{k-r}}{k!}\sum_{j=0}^{r}\binom{r}{j}(-1)^{r-j}H_{n-k+r,\lambda}(j). \nonumber
\end{align}
Therefore, by \eqref{43} and \eqref{46}, we obtain the following theorem.
\begin{theorem}
For $n\ge r$, we have the representation given by
\begin{align*}
H_{n,\lambda}(x)\ &=\sum_{k=0}^{r-1}\frac{1}{k!}\bigg\{\sum_{j=0}^{k}\sum_{l=0}^{n}\frac{2^{l}\binom{n}{l}\binom{k}{j}}{\binom{l+r-k}{l}}S_{2,\lambda}(l+r-k,r-k)(-1)^{k-j}H_{n-l,\lambda}(j)\bigg\}\beta_{k,\lambda}^{(r)}(x) \\
&+n!\sum_{k=r}^{n}\bigg\{\sum_{j=0}^{r}\frac{(-1)^{r-j}\binom{r}{j}2^{k-r}}{k!(n-k+r)!}H_{n-k+r,\lambda}(j)\bigg\}\beta_{k,\lambda}^{(j)}(x).
\end{align*}
\end{theorem}
Let
\begin{equation}
H_{n,\lambda}(x)=\sum_{k=0}^{n}C_{n,k}h_{k,\lambda}^{(r)}(x|u). \label{47}	
\end{equation}
Then, by \eqref{22}, \eqref{23} and \eqref{32}, we have
\begin{align}
C_{n,k}\ &=\ \frac{1}{k!}\bigg\langle \frac{\big(\frac{e_{\lambda}(2t)-u}{1-u}\big)^{r}}{e_{\lambda}\big(\frac{1}{4}(2t)^{2}\big)}(2t)^k\bigg|(x)_{n,\lambda}\bigg\rangle_{\lambda}\label{48}\\
&=\ 	\frac{1}{k!}\bigg\langle \frac{\big(\frac{e_{\lambda}(t)-u}{1-u}\big)^{r}}{e_{\lambda}\big(\frac{1}{4}t^{2}\big)}t^{k}\bigg|2^{n}(x)_{n,\lambda}\bigg\rangle_{\lambda}\nonumber \\
&=\ \frac{1}{k!}\bigg\langle \bigg(\frac{e_{\lambda}(t)-u}{1-u}\bigg)^{r}t^{k}\bigg|\bigg(e_{\lambda}^{-1}\bigg(\frac{1}{4}t^{2}\bigg)\bigg)_{\lambda}2^{n}(x)_{n,\lambda}\bigg\rangle_{\lambda}\nonumber \\
&=\frac{1}{k!(1-u)^{r}}\big\langle\big(e_{\lambda}(t)-u\big)^{r}\big|(t^{k})_{\lambda}H_{n,\lambda}(x)\big\rangle_{\lambda}\nonumber \\
&=\ \frac{1}{(1-u)^{r}k!}\big\langle\big(e_{\lambda}(t)-u\big)^{r}\big|2^{k}(n)_{k}H_{n-k,\lambda}(x)\big\rangle_{\lambda}\nonumber \\
&=\frac{\binom{n}{k}2^{k}}{(1-u)^{r}}\sum_{j=0}^{r}\binom{r}{j}(-u)^{r-j}\big\langle e_{\lambda}^{j}(t)\big|H_{n-k,\lambda}(x)\big\rangle_{\lambda} \nonumber \\
&=\frac{\binom{n}{k}2^{k}}{(1-u)^{r}}\sum_{j=0}^{r}\binom{r}{j}(-u)^{r-j}H_{n-k,\lambda}(j).\nonumber
\end{align}
Therefore, by \eqref{47} and \eqref{48}, we obtain the following theorem.
\begin{theorem}
For $n\ge 0$, we have the representation given by
\begin{displaymath}
H_{n,\lambda}(x)=\frac{1}{(1-u)^{r}}\sum_{k=0}^{n}\binom{n}{k}2^{k}\bigg[\sum_{j=0}^{r}\binom{r}{j}(-u)^{r-j}H_{n-k,\lambda}(j)\bigg]h_{k,\lambda}^{(r)}(x|u).
\end{displaymath}
\end{theorem}

Finally, we would like to obtain another expression of $H_{n,\lambda}(x)$ in terms of $h_{k,\lambda}^{(r)}(x|u)$. With $C_{n,k}$ as in \eqref{47}, we have
\begin{align*}
C_{n,k}\ &=\ \frac{1}{k!}\bigg\langle\frac{\big(\frac{e_{\lambda}(2t)-u}{1-u}\big)^{r}}{e_{\lambda}\big(\frac{1}{4}(2t)^{2}\big)}(2t)^{k}\big|(x)_{n,\lambda}\big\rangle_{\lambda} \\
&=\frac{2^{k}(n)_{k}}{k!(1-u)^{r}}\bigg\langle\frac{(e_{\lambda}(2t)-u)^{r}}{e_{\lambda}(t^{2})}\bigg|(x)_{n-k,\lambda}\bigg\rangle_{\lambda}\\
&=\ \frac{\binom{n}{k}2^{k}}{(1-u)^{r}}\sum_{l=0}^{[\frac{n-k}{2}]}\frac{(-1)_{l,\lambda}}{l!}\big\langle \big(e_{\lambda}(2t)-u\big)^{r}\big|(t^{2l})_{\lambda}(x)_{n-k,\lambda}\big\rangle_{\lambda}\\
&=\ \frac{\binom{n}{k}2^{k}}{(1-u)^{r}}\sum_{l=0}^{[\frac{n-k}{2}]}\frac{(-1)_{l,\lambda}(n-k)_{2l}}{l!}\big\langle \big(e_{\lambda}(2t)-u\big)^{r}\big|(x)_{n-k-2l,\lambda}\big\rangle_{\lambda}\\
&=\ \frac{\binom{n}{k}2^{k}}{(1-u)^{r}}\sum_{l=0}^{[\frac{n-k}{2}]}\frac{(-1)_{l,\lambda}(n-k)!}{l!(n-k-2l)!}\sum_{j=0}^{r}\binom{r}{j}(-u)^{r-j}\langle e_{\lambda}^{j}(2t)\big|(x)_{n-k-2l,\lambda}\big\rangle_{\lambda}\\
&=\ \frac{1}{(1-u)^{r}}\sum_{j=0}^{r}\sum_{l=0}^{[\frac{n-k}{2}]}\frac{\binom{n}{k}\binom{r}{j}2^{k}(-1)_{l,\lambda}(n-k)!}{l!(n-k-2l)!}(-u)^{r-j}2^{n-k-2l}(j)_{n-k-2l,\lambda}
\end{align*}

Hence, we have the following theorem.
\begin{theorem}
For $n\ge 0$, we have the representation given by
\begin{displaymath}
H_{n,\lambda}(x)=\frac{1}{(1-u)^{r}}\sum_{k=0}^{n}\bigg[\sum_{j=0}^{r}\sum_{l=0}^{[\frac{n-k}{2}]}\frac{\binom{n}{k}\binom{r}{j}(-1)_{l,\lambda}(n-k)!}{l!(n-k-2l)!}(-u)^{r-j}2^{n-2l}(j)_{n-k-2l,\lambda}\bigg]h_{k,\lambda}^{(r)}(x).
\end{displaymath}
\end{theorem}

\section{conclusion}

Special polynomials and numbers can be studied by several different methods, which include generating functions, combinatorial methods, umbral calculus, $p$-adic analysis, differential equations, probability, orthogonal polynomials and special functions. These various means of investigating special polynomials and numbers can be applied also to degenerate special polynomials and numbers. Indeed, in recent years, degenerate versions of many special polynomials and numbers were explored with such methods, and some of their arithmetical and combinatorial properties were discovered. Moreover, those degenerate versions of some special polynomials found some applications to other areas of mathematics such as differential equations, identities of symmetry and probability theory.  \par
Hermite polynomials are orthogonal polynomials that arise in such diverse areas as combinatorics, numerical analysis, probability, physics, random matrix theory and systems theory. In light of the regained recent interests in degenerate special numbers and polynomials, the introduction of the degenerate Hermite polynomials followed naturally. In [8], the $\lambda$--linear functionals and the $\lambda$--differential operators were introduced in order to effectively treat $\lambda$--Sheffer polynomials (see \eqref{17}, \eqref{18}). In particular, the important formula in \eqref{22} and \eqref{23}, which expresses one $\lambda$--Sheffer polynomial in terms of other $\lambda$--Sheffer polynomials, were derived in [8]. In this paper, we applied this formula and expressed the higher-order degenerate Bernoulli, Euler and Frobenius-Euler polynomials respectively in terms of
the degenerate Hermite polynomials and vice versa. \par
As one of our future research projects, we would like to continue to study `$\lambda$-umbral calculus' and
their applications to physics, science and engineering as well as to mathematics.

\vspace{0.2cm}

\noindent{\bf{Acknowledgments}} \\
The authors would like to thank Jangjeon Research Institute for Mathematical Sciences for its support during the preparation of this paper.\\

\vspace{0.05cm}

\noindent{\bf{Funding}} \\
Not applicable. \\

\vspace{0.05cm}

\noindent{\bf{Ethics approval and consent to participate}}\\
All authors reveal that there is no ethical problem in the production of this paper. \\

\vspace{0.05cm}

\noindent{\bf{Competing interests}} \\
The authors declare that they have no competing interests. \\

\vspace{0.05cm}

\noindent{\bf{Consent for publication}} \\
All authors want to publish this paper in this journal. \\

\vspace{0.05cm}

\noindent{\bf{Authors’ contributions}} \\
TK and DSK conceived of the framework and structured the whole paper; TK wrote the paper; L-CJ, HL and HYK checked the results of the paper; DSK completed the revision of the article. All authors read and approved the final manuscript. \\

\vspace{0.05cm}

\noindent{\bf{Author details}}\\


\begin{thebibliography}{9}
\bibitem{1}
Araci, S. \emph{Novel identities involving Genocchi numbers and polynomials arising from applications of umbral calculus,} Appl. Math. Comput. \textbf{233} (2014), 599–-607.
\bibitem{2}
Bell, E. T. \emph{Umbral symmetric functions and algebraic analogues of the Bernoulli an and Eulerian numbers and functions,} Math. Z. \textbf{19} (1924), no. 1, 35–-49.
\bibitem{3}
Carlitz, L. \emph{Degenerate Stirling, Bernoulli and Eulerian numbers,} Utilitas Math. \textbf{15} (1979), 51-–88.
\bibitem{4}
Carlitz, L. \emph{A degenerate Staudt-Clausen theorem,} Arch. Math. (Basel) \textbf{7} (1956), 28-–33.
\bibitem{5}
Cigler, J. \emph{Some remarks on Rota's umbral calculus,} Nederl. Akad. Wetensch. Proc. Ser. A 81=Indag. Math. \textbf{40} (1978), no. 1, 27-–42.
\bibitem{6}
Dere, R.; Simsek, Y. \emph{Hermite base Bernoulli type polynomials on the umbral algebra,} Russ. J. Math. Phys. \textbf{22} (2015), no. 1, 1-–5.
\bibitem{7}
Dere, R.; Simsek, Y. \emph{Applications of umbral algebra to some special polynomials,} Adv. Stud. Contemp. Math. (Kyungshang) \textbf{22} (2012), no. 3, 433-–438.
\bibitem{8}
Kim, D. S.; Kim, T. \emph{Degenerate Sheffer sequences and $\lambda$-Sheffer sequences,} J. Math. Anal. Appl. \textbf{493} (2021), no. 1, 124521.
\bibitem{9}
Kim, D. S.; Kim, T. \emph{A Note on a New Type of Degenerate Bernoulli Numbers,} Russ. J. Math. Phys. \textbf{27} (2020), no. 2, 227-–235.
\bibitem{10}
Kim, T.; Kim, D. S. \emph{Degenerate Laplace Transform and Degenerate Gamma Function,} Russ. J. Math. Phys. \textbf{24} (2017), no. 2, 241--248.
\bibitem{11}
Kim, T.; Kim, D. S. \emph{Degenerate polyexponential functions and degenerate Bell polynomials,} J. Math. Anal. Appl. \textbf{487} (2020), no. 2, 124017, 15 pp.
\bibitem{12}
Kim, T.; Kim, D. S.; Kim, H.-Y.; Lee, H.; Jang, L.-C. \emph{Degenerate Bell polynomials associated with umbral calculus,} J. Inequal. Appl. 2020, 2020:226.
\bibitem{13}
Kim, T.; Kim, D. S.; Lee, H.; Park, J.-W. \emph{A note on degenerate $r$-Stirling numbers,} J. Inequal. Appl. 2020, 2020:225.
\bibitem{14}
Kim, T.; Kwon, H.-I.; Seo, J. J. \emph{Some identities of degenerate Frobenius-Euler polynomials and numbers,} Proc. Jangjeon Math. Soc. \textbf{19} (2016), no. 1, 157–-163.
\bibitem{15}
Ma, Y.; Kim, D. S.; Kim, T.; Kim, H.; Lee, H. \emph{Some identities of Lah-Bell polynomials,} Adv. Difference Equ. 2020, 2020:510.
\bibitem{16}
Ma, Y.; Kim, T. \emph{A note on negative $\lambda$-binomial distribution,} Adv. Difference Equ. 2020, 2020:569.
\bibitem{17}
Parrish, C. \emph{Multivariate umbral-calculus,} Thesis (Ph.D.)–-University of California, San Diego, 1974, 90 pp.
\bibitem{18}
Roman, S. \emph{The umbral calculus,} Pure and Applied Mathematics, 111. Academic Press, Inc. [Harcourt Brace Jovanovich, Publishers], New York, 1984.
\bibitem{19}
Roman, S. \emph{More on the umbral calculus, with emphasis on the $q$-umbral calculus,} J. Math. Anal. Appl. \textbf{107} (1985), no. 1, 222-–254.
\bibitem{20}
Roman, S. The theory of the umbral calculus. III. J. Math. Anal. Appl. \textbf{95} (1983), no. 2, 528-–563.
\end{thebibliography}
\end{document}